\documentstyle[12pt,twoside]{article}

\font\teneufm=eufm10 scaled \magstep1
\font\seveneufm=eufm7 scaled \magstep1
\font\fiveeufm=eufm5  scaled \magstep1
\newfam\eufmfam
\textfont\eufmfam=\teneufm
\scriptfont\eufmfam=\seveneufm
\scriptscriptfont\eufmfam=\fiveeufm
\def\frak#1{{\fam\eufmfam\relax#1}}

\newfam\msbfam
\font\tenmsb=msbm10 scaled \magstep1  \textfont\msbfam=\tenmsb
\font\sevenmsb=msbm7 scaled \magstep1 \scriptfont\msbfam=\sevenmsb
\font\fivemsb=msbm5 scaled \magstep1  \scriptscriptfont\msbfam=\fivemsb
\def\Bbb{\fam\msbfam \tenmsb}

\def\RR{{\Bbb R}}
\def\CC{{\Bbb C}}

\def\NN{{\Bbb N}}
\def\ZZ{{\Bbb Z}}

\def\PP{{\Bbb P}}

\def\ra{\rightarrow}

 \def\HollowBoxx #1#2#3{{\dimen0=#1 \advance\dimen0 by -#2
       \dimen1=#1 \advance\dimen1 by #3
        \vrule height 0pt depth #3 width #2
       \hskip -#3
       \vrule height #1 depth #3 width #3}}
 \def\LeftContraction{\mathord{\kern1.45pt \HollowBoxx{6pt}{3.5pt}{.4pt}}\,}

 \def\HollowBox #1#2#3{{\dimen0=#1 \advance\dimen0 by -#3
       \dimen1=#1 \advance\dimen1 by #3
        \vrule height #1 depth #3 width #3
        \vrule height 0pt depth #3 width #2
        \hskip -#3}}
 \def\RightContraction{\mathord{\, \HollowBox{6pt}{3.1pt}{.4pt}} \kern1.6pt}

\def\qed{{\hfill $\Box$}}
\newtheorem{theorem}{THEOREM}[section]

\newtheorem{lemma}[theorem]{Lemma}
\newtheorem{example}[theorem]{Example}
\newtheorem{remark}[theorem]{Remark}
\newtheorem{proposition}[theorem]{Proposition}

\newtheorem{definition}[theorem]{Definition}

\begin{document}

\begin{center}
{\Large \bf Effective Actions of $SU_n$ on Complex
\medskip\\
$n$-dimensional Manifolds}
\medskip \\
\normalsize A. V. Isaev and N. G. Kruzhilin\footnote{The second author
  is supported by grants RFBR 99-01-00969a and 00-15-96008.}

\end{center}

\begin{quotation} \small \sl For $n\ge 2$ we classify all connected $n$-dimensional
  complex manifolds admitting effective actions of the special unitary
  group $SU_n$ by biholomorphic transformations.
\end{quotation}

\pagestyle{myheadings}
\markboth{A. V. Isaev and N. G. Kruzhilin}{$SU_n$-Actions}

\setcounter{section}{-1}

\section{Introduction}
\setcounter{equation}{0}

Let $M$ be a complex manifold and $\hbox{Aut}\,(M)$ the group of biholomorphic automorphisms of $M$ equipped with the compact-open topology. An action of a Lie group $G$ on $M$ by biholomorphic transformations is a real-analytic map
$$
\Phi: G\times M\ra M,
$$
such that for every $g\in G$ we have $\Phi(g,\cdot)\in \hbox{Aut}\,(M)$, and the induced mapping $\Psi: G\ra\hbox{Aut}\,(M)$, $g\mapsto\Phi(g,\cdot)$ is a homomorphism. In this paper we consider the special case $G=SU_n$.

Actions of the group $SU_n$ on real manifolds have been studied extensively. One motivation for such studies is the importance of $SU_n$-actions in physics, especially for small values of $n$ (see, e. g., \cite{KS}). $SU_n$-actions have also been of interest to mathematicians, and various classification results for such actions have been obtained (see, e.g., \cite{HsiWC}, \cite{HsiWY}, \cite{M}). There is, however,  no classification result for the case of $SU_n$-actions by biholomorphic transformations on complex manifolds. We note, however, that in \cite{U} all real compact connected orientable manifolds of dimension $2n$
admitting actions of $SU_n$ were found for $n\ge 5$.

In the present paper we give a complete classification of complex $n$-dimensional manifolds that admit effective actions of the group $SU_n$ by biholomorphic transformations for $n\ge 2$. The effectiveness of an action means that the map $\Psi$ defined above is injective. For an effective action, $\hbox{Aut}\,(M)$ contains a subgroup isomorphic to $SU_n$.

In Section 1 we consider the simplest case when an action has a fixed point. In this case $M$ is equivalent to either the unit ball $B^n\subset\CC^n$, or $\CC^n$, or $\CC\PP^n$ (Proposition \ref{fixed}). 

The rest of the paper deals with actions without fixed points. In Section 2 we describe orbits of such actions (Theorem \ref{orbitclass}). It turns out that every orbit is either a real, or complex hypersurface in $M$. 

In Section 3 we show how orbits can be glued together. We first consider the case when all orbits are real hypersurfaces and show that for $n\ge 3$
a  manifold that admits such an action is equivalent to either a spherical shell in $\CC^n$, or a Hopf manifold, or the quotient  of one of these manifolds
 by the action of a discrete subgroup of the center of
$U_n$. For $n=2$, however, the situation is more interesting. Apart from the above manifolds the
classification in this case also includes spherical shells in $\CC^2$ with
a non-standard complex structure inherited from
the non-standard complex structure on $\CC\PP^2\setminus\{0\}$
introduced in \cite{R1} (Theorem \ref{finalstep}).

Next, we consider the situation when at least one complex hypersurface orbit is present in $M$ and show that there can exist  at most
 two such orbits. They are
 biholomorphically equivalent to $\CC\PP^{n-1}$ and, for $n\ge 3$, can
only arise  either as a result of blowing up
 $\CC^n$ or a ball in $\CC^n$
 at the origin, or adding the hyperplane
$\infty\in\CC\PP^n$ to the exterior of a ball in $\CC^n$,
or blowing up $\CC\PP^n$ at
one point,  or taking the quotient of any
of these examples by the action  of a discrete subgroup of the center
of $U_n$. For $n=2$ the classification
also includes the exterior of a ball in $\CC\PP^2\setminus\{0\}$
with non-standard complex structure to which
 the hyperplane $\infty\in\CC\PP^2$ is attached (Theorem \ref{complh}).

One can attempt to obtain classifications analogous to ours in more general settings, for example, for the group $SU_n$ acting on $k$-dimensional complex manifolds with $k\ne n$. In fact, it can be shown that effective actions of $SU_n$ do not exist on manifolds of dimension $k<n$. Thus, our classification is obtained for the smallest possible dimension of manifolds for which there are effective actions. Another generalization is possible if one considers not necessarily effective actions, e.g., actions with non-trivial discrete kernel. For many our arguments the effectiveness of actions is essential. For non-effective actions entirely new  effects are possible, for example, a manifold may not have any real hypersurface orbits and, if $n=2$, totally real codimension 2 orbits can occur.

In \cite{IKru} we classified all complex manifolds of dimension $n$ that admit effective actions of the full unitary group $U_n$. Our study in \cite{IKru} was motivated by a characterization of the complex space $\CC^n$ obtained as a result  of the classification. Our original proof for the case of $SU_n$ was similar (however, harder on the technical side) to that for $U_n$. The shorter and more elegant argument presented in this paper is to a great extent due to communications that we have had with A. Huckleberry. He made extensive comments that significantly changed our original approach. We wish to thank him for his interest in our work and many inspiring suggestions.

We acknowledge that a large part of this work was done while the first author
was visiting the University of Illinois at Urbana-Champaign.

\section{Case of Fixed Point}
\setcounter{equation}{0}

In this section we list complex manifolds that admit effective actions of $SU_n$ by biholomorphic transformations with fixed point. As shown in Proposition \ref{fixed} below, the classification in this case easily follows from the results in \cite{GK} and \cite{BDK}. First, we will introduce some notation.

For $p\in M$ let  $I_p$ be  the isotropy subgroup of $SU_n$ at $p$, i.e., $I_p:=\{g\in SU_n: gp=p\}$. As above, we denote by $\Psi$ the continuous homomorphism of $SU_n$ into $\hbox{Aut}(M)$ induced by the action of $SU_n$ on $M$. Let $L_p:=\{d_p(\Psi(g)): g\in I_p\}$ be the linear isotropy subgroup, where $d_pf$ is the differential of a map $f$ at $p$. Clearly, $L_p$ is a compact subgroup of $GL(T_p(M),\CC)$, where $T_p(M)$ is the tangent space to $M$ at $p$. Since the action of $SU_n$ is effective, $L_p$ is isomorphic to $I_p$. The isomorphism is given by the map
\begin{equation}
\delta: I_p\ra L_p,\qquad \delta(g):=d_p(\Psi(g)).\label{linis}
\end{equation}

We will now prove the following proposition.
 
\begin{proposition}\label{fixed} \sl Let $M$ be a
connected complex manifold of dimension $n\ge 2$ endowed with  an effective
action of  $SU_n$ by biholomorphic transformations that has
 a fixed point in $M$. Then $M$ is
  biholomorphically equivalent to
  either

\noindent (i) the unit ball $B^n\subset\CC^n$, or

\noindent (ii) $\CC^n$, or

\noindent (iii) $\CC\PP^n$.

\noindent The biholomorphic equivalence $f$ can be chosen to be either $SU_n$-equivariant, or, if $n\ge 3$, $SU_n$-antiequivariant, i.e., to satisfy either the relation 
\begin{equation}
f(gq)=gf(q),\label{qua}
\end{equation}
or the relation
\begin{equation}
f(gq)=\overline{g}f(q),\label{qua1}
\end{equation}
for all $g\in SU_n$ and $q\in M$ (here manifolds (i)-(iii) are
considered with the standard action of $SU_n$).
\end{proposition}

\noindent {\bf Proof:} Let $p\in M$ be a fixed point for the $SU_n$-action. Then $I_p=SU_n$, and 
$L_p$ is a subgroup of $GL(T_p(M),\CC)$ isomorphic to $SU_n$. Since $L_p$ is compact, one can find coordinates
in $T_p(M)$ such that $L_p\subset U_n$.  In these coordinates $L_p=SU_n$ and therefore $L_p$ acts transitively on the unit sphere in $T_p(M)$.

Assume first that $M$ is non-compact. Then by \cite{GK} the manifold $M$ is
biholomorphically equivalent to either $B^n$, or $\CC^n$ and a
biholomorphism may be chosen to satisfy $F(gq)=\gamma(g)F(q)$ for
all $g\in SU_n$ and $q\in M$,
and some automorphism $\gamma$ of $SU_n$, where  the
 action of $SU_n$ on $\CC^n$
in the right-hand side is  standard.
Every automorphism of $SU_n$ has either the form
\begin{equation}
g\mapsto h_0gh_0^{-1},\label{autoform1}
\end{equation}
or the form
\begin{equation}
g\mapsto h_0\overline{g}h_0^{-1},\label{autoform2}
\end{equation}
with $h_0\in SU_n$.
Thus, setting $f:=\hat h_0^{-1}\circ F$, where $\hat h_0$ is the automorphism of $\CC\PP^n$ corresponding to $h_0$,
we obtain either (\ref{qua}), or (\ref{qua1}) respectively.

Assume  now that $M$ is compact. Then by \cite{BDK} $M$ is
biholomorphically equivalent to $\CC\PP^n$. As above, we denote the equivalence map by $F$. We will now show that a biholomorphism between $M$ and $\CC\PP^n$ can be chosen to satisfy (\ref{qua}) or (\ref{qua1}).

The action of $SU_n$ on $M$ induces an injective homomorphism
$\tilde\Psi:SU_n\ra\hbox{Aut}(\CC\PP^n)$. Since $\tilde\Psi(SU_n)$ has a fixed point in
$\CC\PP^n$, $\tilde\Psi(SU_n)$ is conjugate in $\hbox{Aut}(\CC\PP^n)$ to
$SU_n$ embedded in $\hbox{Aut}(\CC\PP^n)$ in the standard way. Hence there
exists an automorphism $\gamma$ of $SU_n$ such that
 for
some $s\in \hbox{Aut}(\CC\PP^n)$ we have
$(s\circ F)(gq)=\gamma(g)(s\circ F)(q)$ for all $g\in SU_n$ and $q\in M$,  where
 the  action of $SU_n$ on $\CC\PP^n$ in the
right-hand side is standard.
We again use that
 $\gamma$ has an  explicit expression
as in (\ref{autoform1}) or  (\ref{autoform2}) and setting
$f:=\hat h_0^{-1}\circ s\circ F$ obtain a map that satisfies either
(\ref{qua}), or (\ref{qua1}) respectively.

The proof is complete. \qed
\smallskip\\

\section{Description of Orbits}
\setcounter{equation}{0}

In this section we assume that $SU_n$ acts on $M$ without fixed points and give a description of orbits that such $SU_n$-actions can have. We start with examples.
\smallskip\\

\begin{example}\label{threetypes}\rm\hfill

(I). Denote by $\CC^n\setminus\{0\}/\ZZ_m$, with $m\in\NN$, the manifold obtained from $\CC^n\setminus\{0\}$  by identifying every point $z$ in $\CC^n\setminus\{0\}$ with $e^{\frac{2\pi i}{m}}z$. 
Let $SU_n$ act on  $\CC^n\setminus\{0\}/\ZZ_m$ in the standard way. Then the lens manifold ${\cal L}^{2n-1}_m:=S^{2n-1}/\ZZ_m$ is an orbit of this action.

(II). We recall the example of a non-standard complex structure on $\CC\PP^2\setminus\{0\}$ given by Rossi in \cite{R1}. Let $(w_0:w_1:w_2:w_3)$ be homogeneous coordinates in $\CC\PP^3$. Consider in $\CC\PP^3$ the variety ${\cal W}$ given by
\begin{equation}
w_1w_2=w_3(w_3+w_0).\label{W}
\end{equation}
Let $(z_0:z_1:z_2)$ denote homogeneous coordinates in $\CC\PP^2$. Consider the map $\pi: \CC\PP^2\setminus\{0\}\ra {\cal W}$ defined by the formulas
\begin{equation}
\begin{array}{l}
w_0=z_0^2,\\
\displaystyle w_1=z_1^2-\frac{z_1\overline{z_2}}{|z_1|^2+|z_2|^2}z_0^2,\\
\displaystyle w_2=z_2^2+\frac{\overline{z_1}z_2}{|z_1|^2+|z_2|^2}z_0^2,\\
\displaystyle w_3=z_1z_2-\frac{|z_2|^2}{|z_1|^2+|z_2|^2}z_0^2.
\end{array}\label{pi}
\end{equation}
The map $\pi$ is everywhere 2-to-1, and its image is ${\cal W}\setminus \Gamma$, where $\Gamma$ is given by
\begin{equation}
w_0=1,\qquad, w_2=-\overline{w_1},\qquad w_3\in\RR.\label{Gamma}
\end{equation}

Consider the unique complex structure on $\CC\PP^2\setminus\{0\}$ that makes $\pi$ locally biholomorphic. Denote $\CC\PP^2\setminus\{0\}$ with this new complex structure by ${\cal X}$. It can be checked that the standard action of $SU_2$ on ${\cal X}$  is in fact an action by biholomorphic transformations. Denote by ${\frak S}^3_R$ the sphere of radius $R$ in ${\cal X}$ with the induced CR-structure. It is an orbit under the action of $SU_2$ on ${\cal X}$ and therefore its CR-structure is invariant under the standard action of $SU_2$ on the sphere. It follows from the results in \cite{R1} (see also \cite{R2}) that none of the surfaces ${\frak S}^3_R$ is CR-equivalent to the ordinary sphere $S^3$ and hence none of ${\frak S}^3_R$ is spherical, unlike the orbit in (I) above. Further, it can be shown (for example, by using the approach that utilizes classifying algebras as in \cite{Kr}) that every CR-structure on $S^3$ invariant under a transitive action of $SU_2$ by CR-transformations is equivalent to either $S^3$ equipped with the standard CR-structure, or to ${\frak S}^3_R$ for some $R>0$ by means of an $SU_2$-equivariant CR-diffeomorphism, and that the manifolds ${\frak S}^3_R$ are not pairwise CR-equivalent.

(III). Let  $\widehat{\CC^n}$ be the blow-up of $\CC^n$ at the origin, i.e.,
$$
\widehat{\CC^n}:=\left\{(z,w)\in \CC^n\times\CC\PP^{n-1}:z_iw_j=z_jw_i,\,\,\hbox{for
      all $i,j$}\right\},
$$
where $z=(z_1,\dots,z_n)$ are coordinates in $\CC^n$ and
$w=(w_1:\dots:w_n)$ are homogeneous coordinates in
$\CC\PP^{n-1}$. We define an action of $SU_n$ on $\widehat{\CC^n}$ as
follows. For $(z,w)\in\widehat{\CC^n}$ and $g\in SU_n$ we set
$$
g(z,w):=(gz,gw),
$$
where in the right-hand side we use the standard actions of $SU_n$ on $\CC^n$
and $\CC\PP^{n-1}$. Then $\CC\PP^{n-1}$ embedded in $\widehat{\CC^n}$ as the set of all points  $(0,w)\in\widehat{\CC^n}$ is an $SU_n$-orbit.
\end{example}

In this section we will show that every orbit of an $SU_n$-action on $M$ is equivalent to an orbit of one of the tree types specified in Example \ref{threetypes}: a lens manifold ${\cal L}^{2n-1}_m$ for some $m\in\NN$, ${\frak S}^3_R$ for some $R>0$ (here $n=2$) or $\CC\PP^{n-1}$ (see Theorem \ref{orbitclass} for  a precise statement). 

An action of $SU_n$ on $M$ is given by a real-analytic map 
$$
\Phi: SU_n\times M\ra M.
$$
Fix $p\in M$ and let $O(p):=\{gp:g\in SU_n\}$ be the $SU_n$-orbit of $p$. The group $SU_n$ is a totally real submanifold in the complex Lie group $SL_n(\CC)$, and therefore we can locally extend the map $\Phi$ to a holomorphic map
\begin{equation}
\tilde\Phi: V\times M_p\ra M,
\end{equation}
where $V$ is a connected neighborhood of $SU_n$ in $SL_n(\CC)$ and $M_p$ is a neighborhood of $O(p)$ in $M$. We refer to $\tilde\Phi$ as a local holomorphic action of $SL_n(\CC)$ on $M$.  

For a point $p\in M$, let $J_p:=\{g\in V: gp=p\}$  be the local isotropy subgroup of $p$ under the local $SL_n(\CC)$-action. Clearly, $I_p=J_p\cap SU_n$. We now define the normalizer $N_p$ of $J_p$ in $SL_n(\CC)$ as follows (see \cite{Huck}, p. 145): denote by ${\frak j}_p$ the Lie algebra of $J_p$ and set
$$
N_p:=\left\{g\in SL_ n(\CC): g{\frak j}_pg^{-1}={\frak j}_p\right\}.
$$
Clearly, $N_p$ is an algebraic subgroup of $SL_n(\CC)$ and $J_p\subset N_p\cap V$. Further, since we consider actions without fixed points, $N_p$ is a proper subgroup of $SL_n(\CC)$ such that $\hbox{dim}_{\CC}SL_n(\CC)/N_p\le n$. 

We need the following proposition whose proof is similar in part to that of Theorem 3.4 in \cite{Huck} (see p. 169).

\begin{proposition}\label{notfibered} \sl $N_p$ is conjugate in $SL_n(\CC)$ to one of the following subgroups 
\begin{equation}
\left\{\left(
\begin{array}{cc}
1/\det\,A & c\\
0 &\\
\vdots & A\\
0 &
\end{array}
\right),\, A\in GL_{n-1}(\CC),\,c\in\CC^{n-1}\right\},\label{parabolic1}
\end{equation}
\begin{equation}
\left\{\left(
\begin{array}{cccc}
1/\det\,A & 0 & \dots & 0\\
d & & A &
\end{array}
\right),\, A\in GL_{n-1}(\CC),\,d\in\CC^{n-1}\right\}\,\hbox{(here $n\ge 3$)},\label{parabolic2}
\end{equation}
\begin{equation}
S_1:=\left\{
\left(
\begin{array}{cc}
a & 0\\
0 &1/a
\end{array}
\right),
\,\,a\in\CC^{*}
\right\}\,\hbox{(here $n=2$)}, \label{maxred_n=2.1}
\end{equation}
\begin{equation}
S_2:=\left\{
\left(
\begin{array}{cc}
a & 0\\
0 &1/a
\end{array}
\right),
\left(
\begin{array}{cc}
0 & b\\
-1/b & 0
\end{array}
\right),\,\,a,b\in\CC^{*}
\right\}\,\hbox{(here $n=2$)},\label{maxred_n=2.2}
\end{equation}
where in (\ref{parabolic1}) and (\ref{parabolic2}), $c$ is a row vector and $d$ is a column vector respectively.
\end{proposition}

\noindent {\bf Proof:} We say that $SL_n(\CC)/N_p$ cannot be fibered, if there does not exits a proper algebraic subgroup $G\supset N_p$ in $SL_n(\CC)$ such that $\hbox{dim}_{\CC}SL_n(\CC)/G < \hbox{dim}_{\CC}SL_n(\CC)/N_p$ (cf. \cite{Huck}, p. 169). 
Suppose first that $SL_n(\CC)/N_p$ cannot be fibered.  

Assume further that $N_p$ is not reductive and consider its unipotent radical $U$. Since $N_p$ is not reductive, $U$ is non-trivial. Let $N(U)$ be the normalizer of $U$ in $SL_n(\CC)$. Clearly, $N(U)$ is a proper algebraic subgroup of $SL_n(\CC)$. Consider the unipotent radical $W$ of $N(U)$. Suppose first that $W=U$. It then follows from Corollary B in \cite{Hum} (p. 186) that $N(U)$ is parabolic. Since $N(U)\supset N_p$ and $SL_n(\CC)/N_p$ cannot be fibered, we have $\hbox{dim}_{\CC}N(U)=\hbox{dim}_{\CC}N_p$. Since $N(U)$ is connected, we obtain that $N_p=N(U)$ is a maximal proper parabolic subgroup of $SL_n(\CC)$.

Assume now that $U\ne W$. Since $U\subset W$, we have $N_p\subset WN_p$. Further, $WN_p$ is a proper algebraic subgroup of $SL_n(\CC)$ and, since $SL_n(\CC)/N_p$ cannot be fibered, we have $\hbox{dim}_{\CC}WN_p=\hbox{dim}_{\CC}N_p$. Hence $W\subset N_p$ and therefore $W=U$ which is a contradiction. Thus $N_p$ is a maximal proper parabolic subgroup in $SL_n(\CC)$. 

Since $\hbox{dim}_{\CC}SL_n(\CC)/N_p\le n$, $N_p$ is conjugate to either subgroup (\ref{parabolic1}), or, if $n\ge 3$, to subgroup (\ref{parabolic2}), or, if $n=4$, to the subgroup
\begin{equation}
\left\{\left(
\begin{array}{cc}
B & C\\
0 & D
\end{array}
\right),\,B,D\in GL_2(\CC),\, \det B\cdot\det D=1,\, C\in{\frak {gl}}_2(\CC)\right\}.\label{except}
\end{equation}

Let $n=4$ and $N_p$ be conjugate to subgroup (\ref{except}). Since $\hbox{dim}_{\CC}SL_4(\CC)/N_p=4$, we have $\hbox{dim}_{\CC}J_p=\hbox{dim}_{\CC}N_p=11$. Further, $N_p\cap SU_4$ is connected and therefore $I_p=J_p\cap SU_4=N_p\cap SU_4$. Calculating $N_p\cap SU_4$ we obtain that $I_p$ is conjugate in $SU_4$ to $(U(2)\times U(2))\cap SU_4$. It then follows that $SU_4$ acts transitively on $M$. Thus the elements of the center of $SU_4$ act trivially on $M$ and the action in this case  is not effective. Hence $N_p$ for $n=4$ cannot be conjugate to subgroup (\ref{except}).  

Assume next that $N_p$ is reductive and let $K\subset SL_n(\CC)$ be its compact form. Conjugating $K$ if necessary, we can assume that $K\subset SU_n$. Since $\hbox{dim}_{\CC}SL_n(\CC)/N_p\le n$, we have $\hbox{dim}\,K\ge n^2-n-1$. By Lemma 2.1 of \cite{IKra}, such subgroups do not exist for $n\ge 3$. Hence $n=2$, and it follows from Lemma 2.1 of \cite{IKru} that $K^c$, the connected component of the identity of $K$, is conjugate in $SU_2$ to $(U(1)\times U(1))\cap SU_2$. Therefore, $N_p^c$ is conjugate to $S_1$ and thus $N_p$ is conjugate in $SL_2(\CC)$ to either  $S_1$, or $S_2$ (see (\ref{maxred_n=2.1}), (\ref{maxred_n=2.2})). However, $SL_2(\CC)/S_1$ clearly can be fibered since $S_1$ is contained in parabolic subgroup (\ref{parabolic1}). Hence $N_p$ is in fact conjugate to $S_2$.

Suppose now that $SL_n(\CC)/N_p$ can be fibered, i.e., there exists a proper algebraic subgroup $G\supset N_p$ such that $\hbox{dim}_{\CC}SL_n(\CC)/G<\hbox{dim}_{\CC}SL_n(\CC)/N_p\le n$. We can assume that $SL_n(\CC)/G$ cannot be fibered. Arguing as above for $G$ in place of $N_p$ and taking into account that $\hbox{dim}_{\CC}SL_n(\CC)/G<n$, we obtain that $G$ is conjugate in $SL_n(\CC)$ to either subgroup (\ref{parabolic1}), or, if $n\ge 3$, to subgroup (\ref{parabolic2}). In particular, $\hbox{dim}_{\CC}SL_n(\CC)/G=n-1$ and hence $\hbox{dim}_{\CC}J_p=\hbox{dim}_{\CC}N_p=n^2-n-1$. 

Let ${\frak g}$ and ${\frak n}_p$ denote the Lie algebras of $G$ and $N_p$ respectively. Suppose first that $G$ is conjugate to subgroup (\ref{parabolic1}). Then ${\frak g}$ is conjugate in ${\frak {sl}}_n(\CC)$ to the subalgebra 
\begin{equation}
\left\{\left(
\begin{array}{cc}
-\hbox{tr}\, \alpha & \beta\\
0 &\\
\vdots & \alpha\\
0 &
\end{array}
\right),\,\alpha\in{\frak {gl}}_{n-1}(\CC),\, \beta\in\CC^{n-1}\right\},\label{alpha}
\end{equation}
where $\beta$ is a row vector and $\hbox{tr}\,\alpha$ denotes the trace of the matrix $\alpha$. Clearly, ${\frak n}_p$ is a codimension 1 complex subalgebra in ${\frak g}$ which is not an ideal in ${\frak g}$. It is not hard to determine all such subalgebras in ${\frak g}$ to obtain that either $n=2$ and ${\frak n}_p$ is conjugate in ${\frak {sl}}_2(\CC)$ to the subalgebra of diagonal matrices, or $n=3$ and ${\frak n}_p$ is conjugate in ${\frak {sl}}_3(\CC)$ to the subalgebra of upper triangular matrices. We consider these two cases separately.

If $n=2$, $N_p^c$ is conjugate in $SL_2(\CC)$ to subgroup $S_1$. Then $N_p$ is conjugate to either $S_1$, or to $S_2$. Since $SL_2(\CC)/N_p$ can be fibered, $N_p$ is in fact conjugate to $S_1$.

If $n=3$, $N_p^c$ is conjugate in $SL_3(\CC)$ to the subgroup of upper triangular matrices. Therefore $N_p^c\cap SU_3$ is connected and hence $I_p^c=(J_p\cap SU_3)^c=N_p^c\cap SU_3$. Calculating $N_p^c\cap SU_3$ we obtain that $I_p^c$ is conjugate in $SU_3$ to $(U(1)\times U(1)\times U(1))\cap SU_3$. It then follows that $SU_3$ acts transitively on $M$. This implies that the elements of the center of $SU_3$ act trivially on $M$ and the action in this case is not effective. Thus in fact $n\ne 3$.

The case when $G$ is conjugate to subgroup (\ref{parabolic2}) is treated similarly.

The proof is complete. \qed
\smallskip\\

We will now obtain the main result of this section.

\begin{theorem}\label{orbitclass}\sl Let $M$ be a connected complex manifold of dimension $n\ge 2$ endowed with an effective action of $SU_n$ by biholomorphic transformations. Assume that there exist no fixed points for this action. Then for $p\in M$ the orbit $O(p)$ is either a complex, or real hypersurface in $M$. In the first case $O(p)$ is biholomorphically equivalent to $\CC\PP^{n-1}$.
In the second case $O(p)$ is CR-equivalent to either  

\noindent (i) a lens manifold ${\cal L}^{2n-1}_m$ for some $m\in\NN$, $(m,n)=1$, or

\noindent (ii) ${\frak S}^3_R$ for some $R>0$ (here $n=2$).

\noindent The biholomorphic equivalence in the first case and CR-equivalence in the second case can be chosen to be either $SU_n$-equivariant, or, if $n\ge 3$, $SU_n$-antiequivariant (see (\ref{qua}), (\ref{qua1}) respectively).
\end{theorem}

\noindent {\bf Proof:} We apply Proposition \ref{notfibered}. Suppose first that $N_p$ is conjugate to maximal parabolic subgroup (\ref{parabolic1}). Then $SL_n(\CC)/N_p$ is biholomorphically and $SL_n(\CC)$-equivariantly equivalent to $\CC\PP^{n-1}$. Since $\hbox{dim}_{\CC}{\frak {sl}}_n(\CC)/{\frak j}_p\le n$ and $\hbox{dim}_{\CC} SL_n(\CC)/N_p=n-1$, we have either $\hbox{dim}_{\CC}J_p=\hbox{dim}_{\CC}N_p=n^2-n$, or $\hbox{dim}_{\CC}J_p=\hbox{dim}_{\CC}N_p-1=n^2-n-1$.

Let first $\hbox{dim}_{\CC}J_p=n^2-n$. We set
\begin{equation}
\Lambda:=\{gp: g\in V_0\},\label{Lambda}
\end{equation}
\begin{equation}
\Lambda_1:=\{gN_p:g\in V_0\}. \label{Lambda1}
\end{equation}
where $V_0\subset V$ is a connected neighborhood of $SU_n$ in $SL_n(\CC)$. Let $\hat O$ be the $SU_n$-orbit of the element $N_p\in SL_n(\CC)/N_p$. Clearly, $\Lambda$ and $\Lambda_1$ are germs of complex manifolds that contain $O(p)$ and $\hat O$ respectively. 
The map
$$
\rho:\Lambda_1\ra\Lambda,\qquad gN_p\mapsto gp
$$
is well-defined and biholomorphic if $V_0$ is sufficiently small. The restriction of $\rho$ to $\hat O$ is an $SU_n$-equivariant map onto $O(p)$.

Calculating $N_p\cap SU_n$, we see that it is conjugate in $SU_n$ to  the subgroup
\begin{equation}
\left\{\left(
\begin{array}{cc}
1/\det B & 0\\
0 & B
\end{array}
\right),\,B\in U_{n-1}\right\}.\label{diag}
\end{equation}
Since $N_p\cap SU_n$ is connected, we obtain $I_p=J_p\cap SU_n=N_p\cap SU_n$, hence $\hbox{dim}\,O(p)=2n-2$. On the other hand, since $\hbox{dim}_{\CC}J_p=n^2-n$, we have $\hbox{dim}_{\CC}\Lambda=n-1$ and therefore $O(p)=\Lambda$. The same calculation shows that $\hat O=\Lambda_1=SL_n(\CC)/N_p$.  This proves that $O(p)$ is a complex hypersurface in $M$, biholomorphically and $SU_n$-equivariantly equivalent to $\CC\PP^{n-1}$.

Suppose now that $\hbox{dim}_{\CC}J_p=n^2-n-1$. Clearly, ${\frak j}_p$ is a codimension 1 complex normal subalgebra of ${\frak n}_p$. The subalgebra ${\frak n}_p$ is conjugate in ${\frak {sl}}_n(\CC)$ to subalgebra (\ref{alpha}). It is not hard to see that ${\frak j}_p$ is then conjugate in ${\frak {sl}}_n(\CC)$ to the subalgebra of (\ref{alpha}) for which $\alpha\in {\frak {sl}}_{n-1}(\CC)$. Let $H\subset SL_n(\CC)$ be the connected subgroup with Lie algebra ${\frak j}_p$. Clearly, $H$ is conjugate in $SL_n(\CC)$ to the subgroup
$$
\left\{\left(
\begin{array}{cc}
1 & r\\
0 &\\
\vdots & Q\\
0 &
\end{array}
\right),\,Q\in SL_{n-1}(\CC),\,r\in\CC^{n-1}\right\},
$$
where $r$ is a row vector. Further, $H\cap SU_n$ is conjugate in $SU_n$ to $SU_{n-1}$ embedded as the subgroup 
$$
\left\{\left(
\begin{array}{cccc}
1 & 0& \dots & 0\\
0 &&&\\
\vdots & &P&\\
0 &&&
\end{array}
\right),\,P\in SU_{n-1}\right\}.
$$
In particular, $H\cap SU_n$ is connected and therefore for the connected component of the identity $I_p^c$ of $I_p$ we have $I_p^c=H\cap SU_n$. It then follows that $O(p)$ is a real hypersurface in $M$. 

Consider the germ of complex manifold $\Lambda$ defined in (\ref{Lambda}). In this case $\Lambda$ is a neighborhood of $O(p)$ in $M$. Denote now by $\hat O$ the $SU_n$-orbit of the element $H\in SL_n(\CC)/H$ and set
$$
\Lambda_2:=\{gH:g\in V_0\}.
$$
$\hat O$ is a real hypersurface in $SL_n(\CC)/H$ and $\Lambda_2$ is its neighborhood. The holomorphic map
$$
\sigma:\Lambda_2\ra\Lambda,\qquad gH\mapsto gp
$$
is well-defined if $V_0$ is sufficiently small. The restriction $\hat\sigma$ of $\sigma$ to $\hat O$ is an $SU_n$-equivariant covering CR-map onto $O(p)$. The fibers of $\hat\sigma$ are given by $\hat\sigma^{-1}(gp)=gI_pH$ for $g\in SU_n$. It is easy to see that $SL_n(\CC)/H$ is biholomorphically and $SL_n(\CC)$-equivariantly equivalent to $\CC^n\setminus\{0\}$. Hence $\hat O$ is equivalent to the sphere $S^{2n-1}\subset\CC^n$ by means of an $SU_n$-equivariant CR-diffeomorphism. Thus we obtain that $S^{2n-1}$ covers $O(p)$ by means of an $SU_n$-equivariant CR-map $\tilde\sigma: S^{2n-1}\ra O(p)$. 

We will now determine the fibers of $\tilde\sigma$. For this we need to find the full isotropy group $I_p$. Suppose first that $n\ge 3$ and apply Lemma 4.4 of \cite{IKru}. Since $I_p^c$ is conjugate in $SU_n$ to $SU_{n-1}$, we obtain that $I_p$ is conjugate in $SU_n$ to $G^m_n\cdot SU_{n-1}$, $m\in\NN$, where $G^m_n$ is the subgroup
$$
\left\{\left(
\begin{array}{cc}
s & 0\\
0 & t\cdot\hbox{id}
\end{array}
\right),\,s,t\in\CC^*,\, s^m=1,\,st^{n-1}=1\right\}.
$$

We will now show that this also holds for $n=2$. Let $T_p(O(p))$ be the tangent space at $p$ to $O(p)$ in the tangent space $T_p(M)$  at $p$ to $M$. Choose coordinates in $T_p(M)$ in which the linear isotropy subgroup $L_p\subset GL(T_p(M),\CC)$ becomes a subgroup of $U_2$ and consider the  orthogonal complement $W$ to $T_p(O(p))\cap iT_p(O(p))$. Clearly, $\hbox{dim}_{\CC} T_p(O(p))\cap iT_p(O(p))=\hbox{dim}_{\CC}W=1$. The group $L_p$ preserves both $T_p(O(p))\cap iT_p(O(p))$ and $W$. In addition, it preserves $T_p(O(p))$ and hence the line $W\cap T_p(O(p))$. Therefore it can only act as $\pm\hbox{id}$ on $W$. Since $O(p)$ is covered by $S^3$, it is strongly pseudoconvex and therefore $L_p$ can only act trivially on $W$. Thus $L_p$ and hence $I_p$ are isomorphic to a subgroup of $U_1$.  
This implies that $I_p$ is
a finite cyclic group, i.e., $I_p=\{C^l,0\le l< m\}$ for some
$C\in SU_2$ and $m\in \NN$ such that $C^m=\hbox{id}$. Choosing new
coordinates in which $C$ is in the diagonal form  we see that
 $I_p$ is conjugate in $SU_2$ to the group $G^m_2$. 

Thus, we have proved that for all $n\ge 2$, $I_p$ is conjugate in $SU_n$ to $G^m_n\cdot SU_{n-1}$. This implies that the fibers of $\tilde\sigma$ are given as follows: $\tilde\sigma^{-1}(gp)=\{\mu g (1,0,\dots,0): \mu^m=1\}$, where $g (1,0,\dots,0)$ denotes the ordinary action of the element $g\in SU_n$ on the vector $(1,0\dots,0)\in S^{2n-1}$ (here we assume that $\tilde\sigma$ is chosen to satisfy  $\tilde\sigma((1,0,\dots,0))=p$). Since $S^{2n-1}$ covers the lens manifold ${\cal L}^{2n-1}_m$ by means of an $SU_n$-equivariant CR-map with exactly the same fibers, we obtain that $O(p)$ is CR- and $SU_n$-equivariantly equivalent to ${\cal L}^{2n-1}_m$. The $SU_n$-action on $M$ (and hence on $O(p)$) can only be effective, if $(m,n)=1$.

Suppose now that $n\ge 3$ and $N_p$ is conjugate to maximal parabolic subgroup (\ref{parabolic2}). This case is almost identical to the preceding one. The subgroup (\ref{parabolic2}) is mapped into subgroup (\ref{parabolic1}) by the following outer automorphism of $SL_n(\CC)$: $\gamma(g)=(g^{\hbox{\tiny T}})^{-1}$, where $g^{\hbox{\tiny T}}$ denotes the transponed matrix. The restriction of $\gamma$ to $SU_n$ is an outer automorphism of $SU_n$: $\gamma(g)=\overline{g}$. This observation shows that $O(p)$ in the case $\hbox{dim}_{\CC}J_p=n^2-n$ is equivalent to $\CC\PP^{n-1}$ by means of a biholomorphic $SU_n$-antiequivariant map and in the case $\hbox{dim}_{\CC}J_p=n^2-n-1$ to ${\cal L}^{2n-1}_m$ by means of an $SU_n$-antiequivariant CR-diffeomorphism.

Let now $n=2$ and $N_p$ be conjugate to subgroup $S_1$ (see (\ref{maxred_n=2.1})): 
$$
N_p=g_0S_1g_0^{-1},
$$
for some $g_0\in SL_2(\CC)$. Conjugating the above identity by a suitable $t\in SU_2$ and replacing $g_0$ by $g_0s$ for a suitable $s\in S_1$, we obtain
$$
tN_pt^{-1}=h_0S_1h_0^{-1},
$$
where
$$
h_0=
\left(
\begin{array}{cc}
1 & c\\
0 &1
\end{array}
\right),\qquad c\in\RR.
$$

Let first $c\ne 0$. In this case $h_0S_1h_0^{-1}\cap SU_2$ is the center of $SU_2$. Since   $tI_pt^{-1}\subset h_0S_1h_0^{-1}\cap SU_2$, $I_p$ is discrete. Hence $O(p)$ is a real hypersurface in $M$. The effectiveness of the $SU_2$-action on $M$ and hence on $O(p)$ then implies that $I_p$ is in fact trivial, and therefore $O(p)$ is diffeomorphic to $S^3$. Let $\hat O$ be the $SU_2$-orbit of $N_p\in SL_2(\CC)/N_p$. $\hat O$ is a real hypersurface in $SL_2(\CC)/N_p$. Consider $\Lambda$ given by (\ref{Lambda}) and $\Lambda_1$ given by (\ref{Lambda1}). Clearly, $\Lambda_1$ is a neighborhood of $\hat O$ in $SL_2(\CC)/N_p$. Since $\hbox{dim}_{\CC}SL_2(\CC)/N_p=2$, we have $\hbox{dim}_{\CC}J_p=\hbox{dim}_{\CC}N_p=1$ and therefore $\Lambda$ is a neighborhood of $O(p)$ in $M$. The holomorphic map
\begin{equation}
\tau:\Lambda\ra\Lambda_1,\qquad gp\mapsto gN_p\label{tau}
\end{equation}
is well-defined if $V_0$ is sufficiently small. The restriction of $\tau$ to $O(p)$ is a 2-to-1 $SU_2$-equivariant covering CR-map onto $\hat O$. Further, $SL_2(\CC)/N_p$ is equivalent biholomorphically and $SL_2(\CC)$-equivariantly to the quadric ${\cal Q}\subset\CC^3$ given by
\begin{equation}
z_1^2+z_2^2+z_3^2=1\label{Q}
\end{equation}
(see \cite{AHR}). The quadric ${\cal Q}$ is affinely equivalent to the finite part of the quadric ${\cal W}$ defined in (\ref{W}). Therefore $O(p)$ has a non-spherical $SU_2$-invariant CR-structure, and it now follows from the discussion in Example (\ref{threetypes}) (II) that $O(p)$ is equivalent to ${\frak S}^3_R$ for some $R>0$ by means of an $SU_2$-equivariant CR-diffeomorphism.

Assume now that $c=0$. In this case $h_0S_1h_0^{-1}\cap SU_2=S_1\cap SU_2$ is the following subgroup
$$
T_1:=\left\{
\left(
\begin{array}{cc}
\alpha & 0\\
0 & \overline{\alpha}
\end{array}
\right),|\alpha|=1
\right\}.
$$
Therefore, $I_p$ is conjugate in $SU_2$ to $T_1$ and hence $O(p)$ is a 2-dimensional submanifold of $M$. Consider the map $\tau$ defined in (\ref{tau}). Its restriction to $O(p)$ is a diffeomorphism onto $\hat O$. Since the $SU_2$-orbit $\hat O$ of $N_p$ in $SL_2(\CC)/N_p$ is totally real, $O(p)$ is totally real as well.

Let $T_p(O(p))$ denote the tangent space to $O(p)$ at $p$. Since $O(p)$ is totally real, we have $T_p(M)=T_p(O(p))+iT_p(O(p))$. Consider the map $\delta$ defined in (\ref{linis}). The effectiveness of the $SU_2$-action implies that $\delta$ is an isomorphism and therefore $\delta(-\hbox{id})$ is a non-trivial element of the linear isotropy subgroup $L_p$. On the other hand, $-\hbox{id}\in I_p$ and therefore $\delta(-\hbox{id})$ acts trivially on $T_p(O(p))$. Since $\delta(-\hbox{id})$  is a complex linear transformation of $T_p(M)$ it is in fact the identity. This contradiction shows that $c\ne 0$ and thus $O(p)$ is CR- and $SU_2$-equivariantly equivalent to ${\frak S}^3_R$ for some $R>0$.
 
Let now $n=2$ and $N_p$ be conjugate to subgroup $S_2$ (see (\ref{maxred_n=2.2})). This case is treated similarly to the preceding one. If $c\ne 0$, $h_0S_2h_0^{-1}\cap SU_2$  
is isomorphic to $\ZZ_4$ and consists of the following four elements
$$
\left\{\left(
\begin{array}{cc}
\pm 1 & 0\\
0 & \pm 1
\end{array}
\right),
\left(
\begin{array}{cc}
\pm ic/\sqrt{1+c^2} & \pm i/\sqrt{1+c^2}\\
\pm i/\sqrt{1+c^2} & \mp ic/\sqrt{1+c^2}
\end{array}
\right)\right\}.
$$
It then follows that $I_p$ is discrete. Hence $O(p)$ is a real hypersurface in $M$. All non-trivial subgroups of $h_0S_2h_0^{-1}\cap SU_2$ contain the center of $SU_2$. The effectiveness of the $SU_2$-action on $M$ and hence on $O(p)$ then implies that $I_p$ is in fact trivial and therefore $O(p)$ is diffeomorphic to $S^3$. Consider the map $\tau$ (see (\ref{tau})). The restriction of $\tau$ to $O(p)$ is a 4-to-1 $SU_2$-equivariant covering CR-map onto $\hat O$. Further, $SL_2(\CC)/N_p$ is biholomorphically and $SL_2(\CC)$-equivariantly equivalent to ${\cal Q}/\ZZ_2$ which as before implies that $O(p)$ is equivalent to ${\frak S}^3_R$ for some $R>0$ by means of an $SU_2$-equivariant CR-diffeomorphism.

If $c=0$, $h_0S_2h_0^{-1}\cap SU_2=S_2\cap SU_2$ is the following subgroup
$$
T_2:=\left\{
\left(
\begin{array}{cc}
\alpha & 0\\
0 & \overline{\alpha}
\end{array}
\right),
\left(
\begin{array}{cc}
0 & \beta\\
-\overline{\beta} & 0
\end{array}
\right), |\alpha|=|\beta|=1
\right\}.
$$
Therefore, $I_p$ is conjugate in $SU_2$ to either $T_1$, or to $T_2$. Hence $O(p)$ is a 2-dimensional submanifold of $M$. The restriction of $\tau$ to $O(p)$ is a map onto $\hat O$ which is 2-to-1 if $I_p$ is conjugate to $T_1$ and a diffeomorphism if $I_p$ is conjugate to $T_2$. As before, this gives that $O(p)$ is totally real which contradicts the effectiveness of the $SU_2$-action. Hence in fact $c\ne 0$.

The proof is complete.\qed
\smallskip\\

\section{Classification of actions without fixed points}
\setcounter{equation}{0}

In this section we obtain a complete classification of connected $n$-dimensional complex manifolds that admit effective actions of $SU_n$ by biholomorphic transformations. As in the preceding section, we consider actions without fixed points.

We start with the case when all orbits are real hypersurfaces. 

\begin{definition}\label{realhyp}\sl Let $S_{r,R}^n:=\{z\in\CC^n:r<|z|<R\}$, $0\le r<R\le\infty$, be
a spherical shell in $\CC^n$. Further, we denote by ${\frak S}^2_{r,R}$,  $0\le r<R\le\infty$,
 the spherical shell $S^2_{r,R}$ equipped with the non-standard complex
 structure  induced by the complex structure of ${\cal X}$ (see Example \ref{threetypes}). Finally, 
for $d\in\CC^*$, $|d|\ne 1$, denote by $M_d^n$ the Hopf
manifold obtained by identifying $z\in\CC^n\setminus\{0\}$ with $d\cdot z$.
\end{definition}

We will now prove the following theorem.

\begin{theorem}\label{finalstep} \sl Let $M$ be a connected complex manifold
of dimension $n\ge 2$ endowed  with an  effective action of  $SU_n$
   by biholomorphic transformations. Assume that
  all orbits of this action are real
  hypersurfaces. Then $M$ is biholomorphically equivalent to either

\noindent (i) $S_{r,R}^n/\ZZ_m$, or

\noindent (ii) $M_d^n/\ZZ_m$, or

\noindent (iii) ${\frak S}^2_{r,R}$ (here $n=2$),

\noindent for some $0\le r<R\le\infty$, $d\in\CC^*$, $|d|\ne 1$, $m\in\NN$, $(m,n)=1$. The biholomorphic equivalence can be chosen to be either $SU_n$-equivariant, or, if $n\ge 3$, $SU_n$-antiequivariant (here  manifolds (i)-(iii) are considered with the standard $SU_n$-actions).    
\end{theorem}

\noindent {\bf Proof:} If $n\ge 3$ or $n=2$ and there exists a spherical orbit, i.e., an orbit equivalent to a lens manifold, then, repeating the proof of Theorem 2.7 in \cite{IKru}, we obtain that $M$ is biholomorphically equivalent to either $S_{r,R}^n/\ZZ_m$, or $M_d^n/\ZZ_m$ by means of an $SU_n$-equivariant or $SU_n$-antiequivariant map.

Suppose now that $n=2$ and the orbit of  every point in $M$ is  non-spherical.
Assume first that $M$ is non-compact. Let $p\in M$.
Then there exists $\rho>0$ such that $O(p)$
is equivalent to ${\frak S}^3_{\rho}\subset {\cal X}$ by means of an $SU_2$-equivariant
CR-diffeomorphism $f$. The map $f$ extends to a biholomorphic
$SU_2$-equivariant map between a neighborhood $U$ of $O(p)$ ($U$ can be
taken to be a connected union of orbits) and
${\frak S}^2_{\rho_1,\rho_2}\subset {\cal X}$ with $0\le\rho_1<\rho<\rho_2\le\infty$.

Let $D$ be a maximal domain in $M$ such that there exists an
$SU_2$-equivariant biholomorphic map $f$ from $D$ onto a spherical shell in ${\cal X}$. Let this shell be ${\frak S}^2_{\rho',\rho''}$ for some $0\le\rho'<\rho''\le\infty$. As shown  above,
such a domain $D$ exists. Assume that $D\ne M$ and let $x$ be a
boundary point of $D$. Consider the orbit $O(x)$. Since $O(x)$ is
non-spherical, there exists an $SU_2$-equivariant CR-diffeomorphism $h$ from $O(x)$
onto ${\frak S}^3_{\tilde\rho}$ for some $\tilde\rho>0$.
This  diffeomorphism extends to an $SU_2$-equivariant biholomorphic map between
a neighborhood $V$ of $O(x)$ (that can be taken to be a union of orbits) and
${\frak S}^2_{\tilde\rho_1,\tilde\rho_2}$ for some
$0\le\tilde\rho_1<\tilde\rho<\tilde\rho_2\le\infty$. For  $s\in V\cap D$ we  consider the orbit $O(s)$. The CR-diffeomorphisms $f$ and $h$ map $O(s)$ into some surfaces ${\frak S}^3_{r_1}$ and ${\frak S}^3_{r_2}$.
Hence the CR-diffeomorphism $F:=h\circ f^{-1}$ maps ${\frak S}^3_{r_1}$ $SU_2$-equivariantly
onto ${\frak S}^3_{r_2}$. Since the surfaces ${\frak S}^3_{r_1}$ and ${\frak S}^3_{r_2}$ are not CR-equivalent unless $r_1=r_2$, it follows that $r_1=r_2=t$, and $F$ is an $SU_n$-equivariant holomorphic automorphism of ${\frak S}^3_t$.

We now need the following lemma.

\begin{lemma}\label{automgp}\sl For any $t>0$, every holomorphic automorphism of ${\frak S}^3_t$ extends to an automorphism of the finite part ${\cal X}'$ of ${\cal X}$, namely ${\cal X}':={\cal X}\setminus\{(0:z_1:z_2)\}$. 
\end{lemma}

\noindent {\bf Proof:} Fix $p\in {\frak S}^3_t$. Since ${\frak S}^3_t$ is real-analytic and strongly pseudoconvex, there
exist local coordinates $(z,w=u+iv)$ near $p$ in which
the equation of ${\frak S}^3_t$ is given in the Chern-Moser normal form
\cite{CM}:
\begin{equation}
v=|z|^2+\sum_{k\ge 2, l\ge 2}F_{k\overline{l}}(z,\overline{z}, u),\label{normal}
\end{equation}
where $F_{k\overline{l}}$ denote terms of order $k$ in $z$ and order
$l$ in $\overline{z}$, and the following normalization conditions
hold
$$
F_{2\overline{2}}\equiv 0,\qquad
F_{2\overline{3}}\equiv 0,\qquad
F_{3\overline{3}}\equiv 0.
$$
Since  ${\frak S}^3_t$ is homogeneous and not spherical, $F_{2\overline{4}}\not\equiv 0$.

Consider the Lie group $\hbox{Aut}({\frak S}^3_t)$ of all 
holomorphic automorphisms of ${\frak S}^3_t$ and denote by $\hbox{Aut}_p({\frak S}^3_t)$ the isotropy subgroup of $p$. Since  ${\frak S}^3_t$ is not spherical at $p$, by \cite{KL} in some normal coordinates near $p$ all elements of
$\hbox{Aut}_p(O(p))$ can be written in the form
\begin{equation}
z\mapsto e^{i\alpha}z,\qquad w\mapsto w,\label{uni}
\end{equation}
where $\alpha\in\RR$. Observe now that among all transformations of the form (\ref{uni}), equation (\ref{normal}) with $F_{2\overline{4}}\not\equiv 0$ can only be invariant under
$$
z\mapsto \pm z,\qquad w\mapsto w.
$$  

Thus for every $p$ the isotropy subgroup $\hbox{Aut}_p({\frak S}^3_t)$ consists of no more than two elements. Since $\hbox{Aut}({\frak S}^3_t)$ is transitive on ${\frak S}^3_t$, it follows that  $\hbox{Aut}({\frak S}^3_t)$ has either one, or two connected components respectively. Let $G_0$ denote the connected component of the identity. Since $\hbox{Aut}_p({\frak S}^3_t)$ for every $p$ is discrete, we have $\hbox{dim}\,\hbox{Aut}({\frak S}^3_t)=3$, and hence $G_0$ consists exactly of the automorphisms induced by the standard action of $SU_2$ on ${\cal X}$. We will now show that $\hbox{Aut}({\frak S}^3_t)$ has indeed another connected component (that we denote by $G_1$) and describe it. We will find a holomorphic automorphism $f$ of ${\frak S}^3_t$ such that $f\not\in G_0$. Then $G_1=\{fg:g\in G_0\}$. 

Recall that $\pi$ defined in (\ref{pi}) is a 2-to-1 covering map from ${\cal X}$ onto ${\cal W}\setminus\Gamma$, where ${\cal W}$ a quadric defined in (\ref{W}) and $\Gamma$ is an exceptional set defined in (\ref{Gamma}). The restriction of $\pi$ to ${\cal X}'$ is a covering map onto ${\cal W}'\setminus\Gamma$, where ${\cal W}'$ is the finite part of ${\cal W}$. It is easy to see that ${\cal W}'$ is affinely equivalent to the quadric ${\cal Q}$ introduced in (\ref{Q}), and under the affine equivalence ${\cal W}'\setminus\Gamma$ is mapped onto ${\cal Q}\setminus\RR^3$. Hence there exists a 2-to-1 covering map $\tilde\pi:{\cal X}'\ra {\cal Q}\setminus\RR^3$. It is clear from the definition of $\pi$ that $\tilde\pi(x)=\tilde\pi(y)$ iff $x=\pm y$.

Consider the following automorphism $h$ of ${\cal Q}\setminus\RR^3$:
$$
z_1\mapsto z_1,\quad z_2\mapsto z_2,\quad z_3\mapsto -z_3.
$$
The automorphism $h$ has fixed points in ${\cal Q}\setminus\RR^3$, e.g., $p_0=(\sqrt{2}, i, 0)$. Let $H$ be a lift of $h$ to the universal cover ${\cal X}'$. The map $H$ is an automorphism of ${\cal X}'$ and satisfies $\tilde\pi\circ H=h\circ\tilde\pi$. Let $\tilde\pi^{-1}(p_0)=\{\pm q_0\}$. Then either $H(q_0)=q_0$, or $H(q_0)=-q_0$. In the first case we set $f:=H$, in the second case $f:=-H$. Hence, $f\in\hbox{Aut}\,({\cal X}')$ is non-trivial and has a fixed point in ${\cal X}'$. 

A direct calculation shows that there exists a surjective homomorphism $\phi:SU_2\ra SO_3(\RR)$ such that $\tilde\pi(gq)=\phi(g)\tilde\pi(q)$ for all $q\in {\cal X}'$ and $g\in SU_2$, where $SO_3(\RR)$ acts on $\CC^3$ in the standard way. Therefore, $\tilde\pi$ maps $SU_2$-orbits to $SO_3(\RR)$-orbits. Since $h$ preserves every $SO_3(\RR)$-orbit, $f\in\hbox{Aut}({\frak S}^3_R)$ for every $R>0$. 

It remains to show that $f\not\in G_0$. Indeed, otherwise in a neighborhood of ${\frak S}^3_t$ the automorphism $f$ would coincide with an automorphism induced by an element of $SU_2$, and hence would coincide with it everywhere and thus would not have a fixed point in ${\cal X}'$.

Therefore, $\hbox{Aut}({\frak S}^3_t)$ has indeed two connected components and $G_1=\{fg:g\in G_0\}$. Since $f$ and every $g\in G_0$ extend to automorphisms of ${\cal X}'$, every element of $\hbox{Aut}({\frak S}^3_t)$ does.

The proof is complete.\qed
\smallskip\\

It now follows from Lemma \ref{automgp} that $F$ extends to an automorphism of ${\cal X}'$ (in fact, since $F$ is in addition $SU_2$-equivariant, the proof of Lemma \ref{automgp} implies that $F$ is from the center of $SU_2$ and thus extends to all of ${\cal X}$). Hence
$$
{\cal F}:=\Biggl\{\begin{array}{l}
F\circ f\qquad\,\hbox{on $D$}\\
h\qquad\qquad\hbox{on $V$}
\end{array}
$$
is a holomorphic map
on $D\cup V$, provided that $D\cap V$ is connected. As the proof of Theorem 2.7 in \cite{IKru} shows, $V$ can be chosen so that  $D\cap V$ is
 indeed connected, and  ${\cal F}$ is one-to-one on $D\cup
V$. Hence $D$ is not maximal. This contradiction implies that in fact $D=M$.

Assume now that $M$ is compact and consider  a domain $D$ defined as above.
Since $M$ is compact, $D\ne M$. For  a boundary point $x$ of $D$ we
consider the orbit $O(x)$. Choose a connected neighborhood $V$ of
$O(x)$ as above, and let $V=V_1\cup V_2\cup O(x)$. As in the proof of
Theorem 2.7 in \cite{IKru}, it turns out that $V_j\subset D$, $j=1,2$,
and hence $M=D\cup O(x)$.

We can now extend $f|_{V_1}$ and $f|_{V_2}$ to $SU_2$-equivariant biholomorphic maps
$f_1$ and $f_2$ respectively, that are defined on $V$ and map it onto
spherical shells in ${\cal X}'$. Then  $f_1$ and $f_2$ map $O(x)$ onto
${\frak S}^3_{r_1}$ and ${\frak S}^3_{r_2}$ respectively, for some
$r_1,r_2>0$. Clearly, $r_1\ne r_2$. However,  the surfaces
${\frak S}^3_R$ are not pairwise  CR-equivalent. This contradiction
shows that $M$ cannot be compact.

The proof is complete.\qed
\smallskip\\

We will now turn to the case when at least one complex hypersurface orbit is present. We need the following proposition.

\begin{proposition}\label{two}\sl Let $M$ be a connected complex
manifold of dimension $n\ge 2$ endowed with an effective action of  $SU_n$
by biholomorphic transformations. Suppose that every  orbit is  a
real or complex hypersurface in $M$ and there exists a complex hypersurface orbit. Then there are at most two complex hypersurface orbits.

If, in addition, $n=2$ and there exists a non-spherical real hypersurface orbit in $M$, then there is exactly one complex hypersurface orbit. All sufficiently small tubular neighborhoods of this orbit constructed from an $SU_n$-invariant distance on $M$ are strongly pseudoconcave. 
\end{proposition}

\noindent{\bf Proof:} Suppose first that all orbits in $M$ are complex hypersurfaces. Recall from the proof of Theorem \ref{orbitclass} that the isotropy subgroup of every point in a complex hypersurface orbit is conjugate in $SU_n$ to subgroup (\ref{diag}) and therefore contains the center of $SU_n$. Hence the center of $SU_n$ acts trivially on every complex hypersurface orbit and, since $M$ is a union of such orbits, the $SU_n$-action on $M$ is not effective. This contradiction shows that there is at least one real hypersurface orbit in $M$. It then follows from \cite{N} (see Corollary 5.8 there) that there can exist at most two complex hypersurface orbits in $M$. 

Let $n=2$ and there is a non-spherical real hypersurface orbit in $M$. It then follows from Theorem \ref{finalstep} that in fact every real hypersurface orbit in $M$ is non-spherical. Suppose there exist two complex hypersurface orbits $O_1$ and $O_2$ in $M$. Fix an $SU_2$-invariant distance function on $M$ and consider a tubular $\epsilon$-neighborhood $U_{\epsilon}$ of $O_1$ not containing $O_2$. Clearly, if $\epsilon$ is sufficiently small, $\partial U_{\epsilon}$ is a connected $SU_n$-invariant real hypersurface in $M$ and hence it is a real hypersurface orbit. Therefore, $U_{\epsilon}$ is either strongly pseudoconvex, or strongly pseudoconcave. Suppose that $U_{\epsilon}$ is strongly pseudoconvex. Then blowing
down $O_1$ in $U_{\epsilon}$ we obtain a Stein analytic space with boundary
$\partial U_{\epsilon}$ (see e.g., \cite{GR}). But this is impossible
since it is shown in \cite{R1} (see also \cite{R2}) that none  of
${\frak S}^3_R$ can bound a Stein analytic space. Hence $U_{\epsilon}$
is strongly pseudoconcave. Therefore $M\setminus\overline{U_{\epsilon}}$ is a strongly pseudoconvex neighborhood of $O_2$ which is impossible by the same argument.

The proof is complete.\qed
\smallskip\\    

 To formulate our next result we need the following definition.

\begin{definition}\label{more}\sl Let as before $\widehat{\CC^n}$ denote the blow-up of $\CC^n$ at the origin and, analogously, let $\widehat{B^n}$ and $\widehat{\CC\PP^n}$ denote the blow-ups of the unit ball $B^n$ and $\CC\PP^n$ at the origin respectively. Let further $\widetilde{S_{r,\infty}^n}\subset\CC\PP^n$ for $r>0$ be the union of the spherical shell $S_{r,\infty}^n$ with infinite outer radius and the hypersurface at infinity in $\CC\PP^n$. Similarly, let $\widetilde{{\frak S}^2_{r,\infty}}\subset {\cal X}$ for $r>0$ be the union of
the spherical shell ${\frak S}^2_{r,\infty}\subset {\cal X}$ with infinite outer radius and the hypersurface at infinity in ${\cal X}$.
\end{definition}

We are now ready to formulate our final classification theorem.

\begin{theorem}\label{complh}\sl Let $M$ be a connected complex manifold of dimension $n\ge 2$ endowed with an effective action of  $SU_n$ by biholomorphic transformations. Suppose that each  orbit of this action is either a real, or complex hypersurface and there exists a complex hypersurface orbit. Then $M$ is biholomorphically equivalent to either

\noindent (i) $\widehat{B^n}/\ZZ_m$, or

\noindent (ii) $\widehat{\CC^n}/\ZZ_m$, or

\noindent (iii) $\widehat{\CC\PP^n}/\ZZ_m$, or

\noindent (iv) $\widetilde{S_{r,\infty}^n}/\ZZ_m$, or

\noindent (v) $\widetilde{{\frak S}^2_{r,\infty}}$ (here $n=2$)

\noindent for some $r>0$, $m\in\NN$, $(m,n)=1$. The biholomorphic equivalence can be chosen to be either $SU_n$-equivariant, or, if $n\ge 3$, $SU_n$-antiequivariant (here manifolds (i)-(v) are considered with the standard $SU_n$-actions). 
\end{theorem}

\noindent{\bf Proof:} Assume first that there is only one complex hypersurface orbit $O$. Consider $\tilde M:=M\setminus O$. Since
$\tilde M$ is clearly non-compact, by Theorem \ref{finalstep} the manifold  $\tilde M$ is equivalent to $\hat M$, where
$\hat M$ is either $S_{r,R}^n/\ZZ_m$, or ${\frak S}^2_{r,R}$, for some $0\le r<R\le\infty$, $m\in\NN$, $(m,n)=1$, by means of a biholomorphic map $f:\tilde M\ra\hat M$ that is either $SU_n$-equivariant, or $SU_n$-antiequivariant. We shall  assume that $f$ is $SU_n$-equivariant; the other case can be dealt with in the same manner.

Let $\phi$ and $\psi$ be the ${\frak {su}}_n$-anticanonical maps defined on $M$ and $\hat M$ respectively (see \cite{Huck} p.p. 166-167 for the definition of ${\frak g}$-anticanonical map). They map $M$ and $\hat M$ into a projective space, and are holomorphic on $\tilde M$ and $\hat M$ respectively. {\it A priori} $\phi$ is only a meromorphic map with possible points of indeterminacy in $O$. However, since $O$ is a complex hypersurface in $M$ and $\phi$ is $SU_n$-equivariant, it follows that $\phi$ is in fact holomorphic on all of $M$.

We will be interested in the level sets of $\phi$ and $\psi$ which form $SU_n$-invariant families of analytic subsets in $M$ and $\hat M$ respectively. A direct calculation shows that the level sets of $\psi$ are of the form $L/\ZZ_m$, where $L$ is the intersection of a complex line passing through the origin in $\CC^n$ with either $S_{r,R}^n$, or ${\frak S}^2_{r,R}$  (in the case $\hat M={\frak S}^2_{r,R}$ we set $m=1$). Since $f$ is $SU_n$-equivariant, it maps the intersections of the level sets of $\phi$ with $\tilde M$ into the level sets of $\psi$. Hence the level sets of $\phi$ form an $SU_n$-invariant family of holomorphic curves in $M$, and for every level set $S$, the intersection $S\cap\tilde M$ is biholomorphically equivalent to either an annulus of modulus $(R/r)^m$ (if $0<r<R<\infty$), or a punctured disk (if $r=0$, $R<\infty$ or $r>0$, $R=\infty$), or $\CC^*$ (if $r=0$ and $R=\infty$). On the other hand, $S\cap\tilde M$ is obtained from the holomorphic curve $S$ by deleting the points where $S$ intersects $O$. Hence $S\cap\tilde M$ cannot be equivalent to an annulus which implies that $S$ intersects $O$ at a single point, and we have $r=0$ or $R=\infty$. Let $\{\epsilon_j\}$ be a sequence of positive numbers convergent to 0. For every $j$ we construct, as in the proof of Proposition \ref{two}, a tubular neighborhood $U_{\epsilon_j}$ of $O$. Since $\partial U_{\epsilon_j}$ is an $SU_n$-orbit and    
$f$ is $SU_n$-equivariant, $f(\partial U_{\epsilon_j})$ is either $r_jS^{2n-1}/\ZZ_m$, or ${\frak S}^3_{r_j}$ for some $r_j>0$. As $j\ra\infty$ we have either $r_j\ra 0$, or $r_j\ra\infty$. 

Assume that $r_j\ra 0$ as $j\ra\infty$ (here $r=0$). Then $U_{\epsilon_j}$ is strongly pseudoconvex. It now follows from Proposition \ref{two} that in this case $\hat M\ne{\frak S}^2_{0,R}$. Hence $\hat M=S_{0,R}^n/\ZZ_m$. Let $B_R^n$ be the ball of radius $R$ in $\CC^n$ and $\widehat{B_R^n}$ its blow-up at the origin (see Example \ref{threetypes} (III) for notation).   
Consider the holomorphic embedding
$\nu: S_{0,R}^n/\ZZ_m\ra \widehat{B_R^n}/\ZZ_m$ defined by the formula
$$
\nu([z]):=\{(z,w)\},
$$
where $w=(w_1:\dots:w_n)$ is uniquely determined by the conditions
$z_iw_j=z_jw_i$ for all $i,j$, $[z]\in(\CC^n\setminus\{0\})/\ZZ_m$
is the equivalence class of
the point $z=(z_1,\dots,z_n)\in\CC^n\setminus\{0\}$ and $\{(z,w)\}\in \widehat{B_R^n}/\ZZ_m$ is the equivalence class of the point $(z,w)\in \widehat{B_R^n}$. 
Clearly, $\nu$ is $SU_n$-equivariant.
Now let  $f_{\nu}:=\nu\circ f$. We
claim  that $f_{\nu}$ extends to $O$ to a
biholomorphic $SU_n$-equivariant map of $M$ onto $\widehat{B_R^n}/\ZZ_m$.

Let $\hat O$ be the unique complex hypersurface orbit in $\widehat{B_R^n}/\ZZ_m$. Take $p\in O$ and find the level set $S_p$ of $\phi$ passing through $p$. Let $\hat p$ be the unique point at which $\overline{f_{\nu}(S_p\setminus\{p\})}$ intersects $\hat O$. Define the extension $F_{\nu}$ of $f_{\nu}$ by setting $F_{\nu}(p)=\hat p$. Clearly, $F_{\nu}$ is $SU_n$-equivariant. We must show that it is continuous at every $p\in O$. Let $\{q_j\}$ be a sequence of points in $M$  converging to $p$. Since all accumulation points of the sequence $\{F_{\nu}(q_j)\}$ lie  in $\hat O$ and $\hat O$ is compact, it suffices  to
show that every  convergent subsequence of $\{F_{\nu}(q_j)\}$ converges to $\hat p$. Let a subsequence $\{F_{\nu}(q_{j_k})\}$ converge to $q\in\hat O$. For every $q_{j_k}$ there exists $g_{j_k}\in SU_n$ such that $g_{j_k}q_{j_k}\in S_p$. We select a convergent subsequence $\{g_{j_{k_l}}\}$ and denote its limit by $g$. Then
$\{g_{j_{k_l}}q_{j_{k_l}}\}$ converges to $gp$. Since $gp\in O$ and $g_{j_{k_l}}q_{j_{k_l}}\in S_p$, it follows that $gp=p$, i.e, $g\in I_p$. By definition, $F_{\nu}$ is continuous on $S_p$ and we have $F_{\nu}(g_{j_{k_l}}q_{j_{k_l}})\ra F_{\nu}(p)$. On the other hand, since $F_{\nu}$ is $SU_n$-equivariant, we have $F_{\nu}(g_{j_{k_l}}q_{j_{k_l}})=g_{j_{k_l}}F_{\nu}(q_{j_{k_l}})\ra gq$. Since $g\in I_p$, this implies that $q=F_{\nu}(p)$. Thus $\{F_{\nu}(q_j)\}$ converges to $F_{\nu}(p)$ which shows that $F_{\nu}$ is continuous and therefore holomorphic on $M$. Hence $M$ is equivalent to either  $\widehat{B^n}/\ZZ_m$, or $\widehat{\CC^n}/\ZZ_m$.

Assume that $r_j\ra \infty$ as $j\ra\infty$ (here $R=\infty$). If $\hat M=S_{r,\infty}^n/\ZZ_m$, we consider  the holomorphic embedding
 $\sigma: S_{r,\infty}^n/\ZZ_m\ra \widetilde{S_{r,\infty}^n}/\ZZ_m$
$$
\sigma([z]):=\{(1:z_1:\dots:z_n)\},
$$
where $(z_0:\dots:z_n)$ are homogeneous coordinates in $\CC\PP^n$, the hyperplane at infinity in $\CC\PP^n$ is given by $\{z_0=0\}$, and $\{(1:z_1:\dots:z_n)\}\in\widetilde{S_{r,\infty}^n}/\ZZ_m$ denotes the equivalence class of  
$(1:z_1:\dots:z_n)\in \widetilde{S_{r,\infty}^n}$. By an analogous argument one can now show that  the map $f_{\sigma}:=\sigma\circ f$ extends to $O$ and gives rise to a
biholomorphic map from $M$ onto $\widetilde{S_{r,\infty}^n}/\ZZ_m$. If $\hat M={\frak S}_{r,\infty}^2$
we regard $\sigma$ as a map from ${\frak S}_{r,\infty}^2$ into $\widetilde{{\frak S}_{r,\infty}^2}$ (setting $m=1$) and  the same proof gives that  $f_{\sigma}$ extends to $O$ and establishes biholomorphic equivalence between $M$ and $\widetilde{{\frak S}_{r,\infty}^2}$.

Assume  now that there exist two complex hypersurface orbits
$O_1$ and $O_2$. As above,
we consider the manifold $\tilde M$ obtained from
$M$ by removing $O_1$ and $O_2$.
 For some $m\in\NN$, $(m,n)=1$, and $0\le r<R\le\infty$,  it is
 biholomorphically equivalent
to $S_{r,R}^n/\ZZ_m$ by means of  a map $f$  that is either $SU_n$-equivariant, or $SU_n$-antiequivariant. Arguments very similar to the ones used above show that in this case
we have $r=0$ and $R=\infty$, and the map $f_{\tau}:=\tau\circ f$ extends to a biholomorphic map
$M\ra\widehat{\CC\PP^n}/\ZZ_m$. Here $\tau: (\CC^n\setminus\{0\})/\ZZ_m\ra
\widehat{\CC\PP}^n/\ZZ_m$ is the $SU_n$-equivariant map defined as
follows:
$$
\tau([z]):=\Bigl\{\Bigl((1:z_1:\dots:z_n),w\Bigr)\Bigr\},
$$
where $w=(w_1:\dots:w_n)$ is uniquely determined from the conditions
$z_iw_j=z_jw_i$ for all $i,j$ and $\Bigl\{\Bigl((1:z_1:\dots:z_n),w\Bigr)\Bigr\}\in\widehat{\CC\PP}^n/\ZZ_m$ is the equivalence class of $\Bigl((1:z_1:\dots:z_n),w\Bigr)\in \widehat{\CC\PP}^n$.  

The proof is complete. \qed

\begin{remark}\label{old}\rm In the proof of Theorem \ref{complh} we extended the maps $f_{\nu}$, $f_{\sigma}$ and $f_{\tau}$ along an $SU_n$-invariant family of holomorphic curves. This family was constructed by considering the level sets of the ${\frak {su}}_n$-anticanonical map on $M$. There are two more constructions that lead to the same family of curves. One approach is to use the isotropy subgroups of points in the exceptional orbits as we did in \cite{IKru} for the group $U_n$ (see the proof of Theorem 3.3 there). Another approach comes from \cite{N}: for an $SU_n$-invariant Hermitian metric one considers the collection of all geodesics passing through a fixed point in an exceptional orbit and orthogonal to it.
\end{remark}

{\obeylines
Department of Mathematics
The Australian National University
Canberra, ACT 0200
AUSTRALIA
E-mail: alexander.isaev@maths.anu.edu.au
\hbox{ \ \ }
\hbox{ \ \ }
Department of Complex Analysis
Steklov Mathematical Institute
42 Vavilova St.
Moscow 117966
RUSSIA
E-mail: kruzhil@ns.ras.ru
}

\end{document}